\documentclass[a4paper,11pt]{article} 
\usepackage{float}
\usepackage{caption}
\usepackage{dsfont}
\usepackage{booktabs} 
\usepackage{tikz-qtree}
\usetikzlibrary{arrows.meta}
\usepackage[a4paper, margin=.9in]{geometry}
\usepackage{amsmath,amsthm,amsfonts,amssymb,array}
\usepackage[square,sort,comma,numbers]{natbib}
\usepackage{subeqnarray}
\usepackage{blkarray}
\usepackage[utf8]{inputenc}
\usepackage[english]{babel}
\usepackage{mathabx}
\usepackage{mathrsfs}
\usepackage{colortbl}
\usepackage{graphicx,epsfig}
\newcounter{myeqno}
\usepackage{color}
\usepackage{tcolorbox}
\usepackage{ upgreek }
\usepackage{relsize}
\usepackage{multirow}
\definecolor{shadecolor}{gray}{0.75}

\theoremstyle{definition}

\newtheorem{proposition}{Proposition}[section]

\newtheorem{theorem}{Theorem}[section]
\newtheorem{example}{Example}[section]
\newtheorem{corollary}{Corollary}[section]

\usepackage[misc,geometry]{ifsym}
\date{}
\usepackage[colorlinks=true]{hyperref}
\hypersetup{urlcolor=blue, citecolor=blue}
\newlength{\defbaselineskip}
\newcommand{\setlinespacing}[1]%
{\setlength{\baselineskip}{#1 \defbaselineskip}}

\begin{document}

		\title{\textbf{On the distance spectral gap and construction of \texorpdfstring{$D$}{D}-equienergetic graphs}}

\author{Haritha T$^1$\footnote{harithathottungal96@gmail.com},  Chithra A. V$^1$\footnote{chithra@nitc.ac.in}
	\\ \\ \small 
 1 Department of Mathematics, National Institute of Technology Calicut,\\\small Calicut-673 601, Kerala, India\\ \small}

\maketitle	
\begin{abstract}
Let $D(G)$ denote the distance matrix of a connected graph $G$ with $n$ vertices. The distance spectral gap of a graph $G$ is defined as $\delta_{D^G} = \rho_1 - \rho_2$, where $\rho_1$ and $\rho_2$ represent the largest and second largest eigenvalues of $D(G)$, respectively. For a   $k$-transmission regular graph $G$, the second smallest eigenvalue of the distance Laplacian matrix equals the distance spectral gap of $G$. In this article, we obtain some upper and lower bounds for the distance spectral gap of a graph in terms of the sum of squares of its distance eigenvalues. Additionally, we provide some bounds for the distance eigenvalues and distance energy of graphs. Furthermore, we construct new families of non $D$-cospectral $D$-equienergetic graphs with diameters of $3$ and $4$.
\end{abstract}

\section{Introduction}
Throughout this article we consider simple, connected, and undirected graphs $G$, with vertex set $V(G)= \{v_1, \ldots, v_n\}$ and edge set $E(G)= \{e_1, \ldots, e_m\}.$ The \textit{adjacency matrix} $A(G)$ of a graph $G$ is an $n\times n$ symmetric matrix, with entries $1$ if $v_i$ is adjacent to $v_j$ and $0$ otherwise. The eigenvalues of $A(G)$ ($A$-eigenvalues) are denoted by $\lambda_1\geq \lambda_2\geq\dots\geq \lambda_n$. Let $d_i$ denote the degree of the $i^{th}$ vertex in a graph $G$. 
The \textit{degree matrix} of $G$, denoted as $Deg(G)$ is a diagonal matrix with $d_i$ as the entries on its diagonal.
The \textit{Laplacian matrix} $L(G)$ of $G$ is defined as $L(G)= Deg(G)-A(G).$ An $n\times m$ matrix in which rows and columns are indexed by vertices and edges respectively and entries are $1$ if $v_i$ is incident to $e_j$, equal to $0$ otherwise, is called the \textit{incidence matrix} $Q(G)$ of $G$. The \textit{line graph} of a graph $G$ is a graph with vertex set as the edge set of $G$, where two vertices are adjacent if the corresponding edges of $G$ have an end vertex in common and it is represented by $l(G)$. If $G$ is $r$-regular then the eigenvalues of $A(l(G))$ consist of the simple eigenvalues $2r-2,\, \lambda_j+r-2$, and the eigenvalue $-2$ with multiplicity $m-n$  ($j=2,3,\ldots, n$) \cite{cvet}. The complement of a graph $G$ is denoted by $\overline{G}.$ We represent the \textit{equitable quotient matrix} \cite{brouwer2011spectra} of a partitioned matrix by $M$. We denote by $K_n$ the \textit{complete graph}, $K_{a,b}$  the \textit{complete bipartite graph} with partitions of sizes $a$ and $b$,  $\left(\text{where}\; a+b= n\right)$ and $C_n$ the \textit{cycle graph}, each of size $n$. Throughout, we denote by $J_{n\times m}$ an $n\times m$ matrix with all entries equal to one and by $I_n$ the identity matrix of order $n$. \\
\par The distance between any two vertices $v_i, v_j\in V(G)$, denoted by $d_G(v_i,v_j)$, is defined to be the length of the shortest path connecting them in $G$ \cite{buckley1990distance}. In particular, $d_G(v_i, v_i)=0$ for any $v_i\in V(G)$. The \textit{distance matrix} \cite{graham1971addressing} $D(G)$ of $G$, is the $n\times n$ symmetric matrix with $(i,j)^{th}$ entry being $d_G(v_i, v_j)$. The eigenvalues of $D(G)$ are said to be the distance eigenvalues of $G$ (denoted by $D$-eigenvalues). Let $\rho_1\geq \rho_2 \geq \cdots \geq \rho_n$ be the \textit{$D$-eigenvalues} of $G$. Also, let $|\rho_1^*|\geq |\rho_2^*|\geq \cdots\geq |\rho_n^*|$ be the non increasing arrangement of $D$-eigenvalues. Two non-isomorphic graphs with the same number of vertices are said to be \textit{$D$-cospectral} if they have equal distance spectrum ($D$-spectrum). In \cite{indulal2010distance}, the authors introduced the concept of \textit{distance energy} ($D$-energy), and it is defined as $E_D(G)= \sum_{i=1}^{n} |\rho_i|.$ Two graphs having an equal number of vertices are said to be \textit{D-equienergetic} if they have the same $D$-energy. $D$-cospectral graphs are evidently $D$-equienergetic, but converse need not be true. In general, non $D$-cospectral graphs need not be $D$-equienergetic. Several studies have explored the construction of $D$-equienergetic graphs and $D$-cospectral graphs using various graph operations \cite{MR4755379,central,MR4703142,haritha2}.\\
\par The \textit{transmission} $Tr(v_i)$ of a vertex $v_i$ is the sum of all distances from $v_i$ to all other vertices in $G$, that is, $Tr(v_i) = \sum_{v_j\in V(G)}d_G(v_i, v_j).$ If $Diag(Tr)$ denotes the diagonal matrix of transmissions of vertices in $G$, then the \textit{distance Laplacian matrix} \cite{aouchiche2013two} is defined as $D^L(G)= Diag(Tr)-D(G).$ \\
\par The \textit{subdivision graph} $S(G)$ of $G$ is obtained by subdividing each edge of $G$ exactly once and the \textit{central graph} $C(G)$ of $G$, is obtained from $S(G)$ by joining all the non adjacent vertices in $G$. In \cite{lu2017generalized}, Lu et al. defined
the following: Let $G_1$ and $G_2$ be two vertex disjoint graphs with vertices $n_1,n_2$ and edges $m_1,m_2$, respectively. Then the \textit{subdivision vertex-vertex} join $G_1 \ocirc G_2$ (resp. \textit{subdivision edge-edge} join $G_1 \ominus G_2$) of
$G_1$ and $G_2$ is the graph obtained from $S(G_1)$ and $S(G_2)$ by
joining every vertex in $V(G_1)$ (resp. each
vertex corresponding to the edges of $G_1$) to every vertex in $V(G_2)$ (resp. each vertex corresponding to the edges of $G_2$). Note that the number of vertices of  $G_1 \ocirc G_2$ and $G_1 \ominus G_2$ are the same and  is  $n_1+n_2+m_1+m_2$. In \cite{jahfar2020central} Jahfar et al. defined the following: The \textit{central vertex join} $ G_1\dot{\vee} G_2$ (\textit{central edge join} $ G_1\veebar  G_2$) of $G_1$ and $G_2$ is the graph obtained from $C(G_1)$ and $G_2$ by joining each vertex of $G_1$ (resp. each vertex corresponding to the edges of $G_1$) with every vertex of $G_2$ (resp. every vertex of $G_2$).\\
\par The spectral gap of a graph $G$ is defined as the difference between the two largest eigenvalues of its adjacency matrix. For regular graphs, the spectral gap is equal to the second smallest eigenvalue of the Laplacian matrix, which is known as the algebraic connectivity.
The spectral gap is primarily studied for the class of (connected) regular graphs
as regular graphs with large spectral gaps are known to demonstrate strong connectivity properties, which are significant in several areas of theoretical computer science. For studies on spectral gap see \cite{regspecgap,specgap2}. Inspired by these studies, this article introduces the  concept of the \textit{distance spectral gap} and investigates its properties, including various bounds. \\
\par In Section $2$, we define the spectral gap for the distance matrix of graphs and provide various bounds for the spectral gap of graphs. Section $3$ presents the lower and upper bounds for the distance eigenvalues and the distance energy of the graphs. Section $4$ introduces a new family of non $D$-cospectral $D$-equienergetic graphs with diameter $4$, using the operation subdivision vertex-vertex join between two graphs.
\section{Distance spectral gap of graphs}

The \textit{distance spectral gap} of a graph $G$ is defined as $\delta_{D^G}= \rho_1-\rho_2.$ The distance spectral gap of the complete graph $K_n$ is $n.$ All trees have exactly one positive $D$-eigenvalue, whereas the maximum distance spectral radius is attained in $P_n$ \cite{ruzieh1990distance} and minimum for $K_{1,n-1}$ \cite{yu2012distance}. Hence the distance spectral gap is maximum for $P_n$ ($n\geq 3$) and minimum for $K_{1,n-1}$.

If $G$ is a $k$-transmission regular graph, then the second smallest eigenvalue of $D^L(G)$ is equal to the distance spectral gap of $G.$ In \cite{aouchiche2013two}, the authors proved that the second smallest distance Laplacian eigenvalue of $G$ is at least $n$ with equality if and only if $\overline{G}$ is disconnected. This implies the following result.
\begin{proposition}
    If $G$ is a transmission regular graph with $n$ vertices, then $\delta_{D^G}\geq n$ with equality if and only if $\overline{G}$ is disconnected.
\end{proposition}
\noindent Consider the $D$-eigenvalues of $G$, $\rho_1\geq \cdots \geq \rho_n$, then $\sum_{i=1}^{n}\rho_i= 0.$ Let $S_D(G)= \sum_{i=1}^{n}\rho_i^2= \sum_{i=1}^{n}\sum_{j=1}^{n}d_G(v_i,v_j)^2.$\\

The following theorem presents bounds for distance eigenvalues of a graph $G.$
\begin{theorem}
     For a graph $G$ with $n$ vertices, let $n^{+}$ and $n^{-}$ be the number of positive and negative eigenvalues of $G$, respectively. Then for $1\leq t\leq n,$
     $$-\sqrt{\frac{S_D(G)n^{+}}{(n-t+1)(n-t+1+n^{+})}}\leq \rho_t \leq \sqrt{\frac{S_D(G)n^{-}}{t(t+n^{-})}}.$$
\end{theorem}
\begin{proof} For $1\leq t\leq n$, the right inequality holds when $\rho_t\leq 0.$ Suppose $\rho_t>0$, then from $S_D(G)$, we have\\
    \begin{equation}
        \rho_t^2= S_D(G)-\sum_{\substack{\rho_i>0\\i\neq t }}\rho_i^2-\sum_{\rho_i< 0}\rho_i^2,
    \end{equation}
where $i\neq t$ and $1\leq i\leq n$. For $\rho_i>0$, the RHS of $(1)$ is maximized when $\rho_1= \rho_2= \cdots = \rho_t$ and, for $\rho_i<0$, $\rho_i= -\frac{t\rho_t}{n^{-}}.$\\
Thus, \begin{equation*}
    \begin{aligned}
        \rho_t^2&\leq S_D(G)-(t-1)\rho_t^2-\frac{t^2\rho_t^2}{n^{-}},\\
        \rho_t&\leq \sqrt{\frac{S_D(G)n^{-}}{t(t+n^{-})}}.
    \end{aligned}
\end{equation*}
The left inequality in the theorem holds for $\rho_t\geq 0.$ Let us assume $\rho_t<0,$ in an analogous manner we have,
\begin{equation*}
    \begin{aligned}
        \rho_t^2&= S_D(G)-\sum_{\rho_i>0}\rho_i^2-\sum_{{\substack{\rho_i<0\\i\neq t}}}\rho_i^2\\
        &\leq S_D(G)-(n-t)\rho_t^2-\frac{(n-t+1)^2\rho_t^2}{n^{+}}\\
        &\leq \frac{S_D(G)n^{+}}{(n-t+1)(n-t+1+n^{+})}.
    \end{aligned}
\end{equation*}
Since $\rho_t<0,$ we have \\
$$\rho_t\geq -\sqrt{\frac{S_D(G)n^{+}}{(n-t+1)(n-t+1+n^{+})}}.$$

\end{proof}
As an observation, when $\rho_t>0$ then $n^{-}\leq n-t$ and when $\rho_t<0$ then $n^{+}\leq t-1$ (note that $\rho_1>0$). Then we have
\begin{corollary}\label{specgap}
    Let $G$ be a graph with $n$ vertices, then for $1\leq t\leq n,$
    $$-\sqrt{\frac{S_D(G)(t-1)}{n(n-t+1)}}\leq \rho_t \leq \sqrt{\frac{S_D(G)(n-t)}{nt}}.$$
\end{corollary}
Next proposition provides lower and upper bounds for the distance spectral gap of a graph.
\begin{proposition}\label{specgap1}
    For a graph $G$ with $n$ vertices, 
    $$-\sqrt{\frac{n-2}{2n}S_D(G)}< \delta_{D^G}\leq \sqrt{\frac{n}{n-1}S_D(G)}.$$
\end{proposition}
\begin{proof}
    From Corollary \ref{specgap}, we have\\
    \begin{equation}
        0< \rho_1\leq \sqrt{\frac{S_D(G)(n-1)}{n}}
    \end{equation}
    and
    
    \begin{equation}
      -\sqrt{\frac{S_D(G)}{n(n-1)}}\leq \rho_2 \leq \sqrt{\frac{S_D(G)(n-2)}{2n}}.   
    \end{equation}

    Then from $(2)$ and $(3)$, we get
$$-\sqrt{\frac{n-2}{2n}S_D(G)}< \delta_{D^G}\leq \sqrt{\frac{n}{n-1}S_D(G)}.$$

\end{proof}
In general, if $S_D(G)$ is maximum then $G$ is the path and minimum then $G$ is the complete graph \cite{zhou2007largest}. Also, 

\begin{equation}
    n(n-1)\leq S_D(G) \leq \frac{n^2(n+1)(n-1)}{6}.
\end{equation}
 the following corollary directly follows from Proposition \ref{specgap1} and $(4)$
\begin{corollary}
    For a connected graph $G$ with $n$ vertices, $$-\sqrt{\frac{n(n-2)(n^2-1)}{12}}< \delta_{D^G}\leq n\sqrt{\frac{n(n+1)}{6}}.$$
\end{corollary}

In \cite{MR3610302}, the authors proved that, for a connected graph $G$, $\rho_2\geq -1$, with equality holds if and only if $G\cong K_n.$ This implies the following proposition.
\begin{proposition}\label{rho+1}
    For a connected graph $G$ with $n$ vertices, $\delta_{D^G}\leq \rho_1+1.$
    The equality holds if and only if $G\cong K_n.$
\end{proposition}
From Corollary $1.10$ in \cite{chen2013sharp} and Proposition \ref{rho+1}, we get the following result.
\begin{proposition}
    Let $G$ be a connected graph with $n$ vertices, row sums of $D(G)$ denoted as $D_1\geq D_2\geq \cdots \geq D_n$, and diameter $d,$ then $$\delta_{D^G}\leq \frac{D_i-d+\sqrt{(D_i+d)^2+4d\sum_{l=1}^{i-1}D_l-D_i}+2}{2},$$
    where $1\leq i\leq n$, with equality holds if and only if $G\cong K_n.$
\end{proposition}
From Theorem $2$ in \cite{zhou2010} and Proposition \ref{rho+1}, we get the following result.
\begin{proposition}
    For a graph $G$ with $n$ vertices, minimum degree $d_1$, second minimum degree $d_2$, and diameter $d$, $$\delta_{D^G}\leq \sqrt{\left(dn-\frac{d(d-1)}{2}-1-d_1(d-1)\right)
    \left(dn-\frac{d(d-1)}{2}-1-d_2(d-1)\right)}+1,$$
    with equality if and only if $G\cong K_n.$
\end{proposition}
\section{Lower and upper bounds for distance energy}
This section provides upper and lower bounds for $D$-energy of $G$ in terms of $S_D(G).$

\begin{proposition}
    Let $G$ be a connected graph with $n$ vertices having non zero $D$-eigenvalues and let $|\rho_1^*|\geq |\rho_2^*|\geq \cdots\geq |\rho_n^*|>0$ be the non increasing arrangement of $D$-eigenvalues, where $\rho_1^*= \rho_1.$ Then $$E_D(G)\geq \rho_1+\frac{S_D(G)-\rho_1^2}{|\rho_2^*|}.$$
    The equality holds if and only if $G\cong K_n.$
\end{proposition}
\begin{proof}
     We have $|\rho_1^*|^2+|\rho_2^*|^2+\cdots+|\rho_n^*|^2= S_D(G).$ \\
    Now,  \[|\rho_2^*|\sum_{i=2}^{n}|\rho_i^*|\geq \sum_{i=2}^{n}|\rho_i^*|^2= S_D(G)-\rho_1^{2}.\]
    Therefore, $E_D(G)\geq \rho_1+\frac{S_D(G)-\rho_1^2}{|\rho_2^*|},$ with equality holds if and only if $G\cong K_n.$
\end{proof}

    

\begin{proposition}
For a connected graph $G$ with $n$ vertices and $D$-eigenvalues $\rho_1, \ldots, \rho_n,$ let the first
$p$ of these eigenvalues be positive , where $1\leq p< n$. Then $$E_D(G)\leq \sqrt{nS_D(G)-\frac{4n}{S_D(G)}\left(\rho_1^2+ \cdots  + \rho_p^2 -\frac{S_D(G)}{2}\right)^2}.$$
\end{proposition}
\begin{proof}
Consider \begin{equation*}
        \begin{aligned}
            \rho_1^2+\cdots +\rho_p^2-\frac{S_D(G)}{2}&= \frac{1}{2}\left( \rho_1^2+\cdots +\rho_p^2-(\rho_{p+1}^2+\cdots+ \rho_n^2)\right)\\
            &=\frac{1}{2}\left(\rho_1|\rho_1|+\cdots +\rho_n|\rho_n|\right)
        \end{aligned}
    \end{equation*}
Then the proof follows by similar arguments in Theorem $5$ in \cite{filipovski2021bounds}.

\end{proof}
\section{Distance equienergetic graphs}
This section explores the construction of $D$-equienergetic graphs of diameter $4$.

\begin{theorem}\label{subdivision vertex}
Let $G_i$ be a $r_i$-regular graph with $n_i$ vertices and $m_i$ edges. And let $A_i$ be the adjacency matrix of $G_i$ and  $\{r_i, \lambda_{i2},\ldots, \lambda_{in_{i}}\}$ ($i=1, 2$) be the adjacency spectrum of $G_i$. Then, the distance spectrum of $G_1\ocirc G_2$ consists of the following,

\begin{itemize}
    \item [(i)] $-2(\lambda_{ij}+r_i+1), j=2,3,\ldots,n_i$,
    \item [(ii)] $0$ with multiplicity $m_i-1$,
    \item[(iii)] the four eigenvalues of the equitable quotient matrix of $D(G_1\ocirc G_2)$,
\[N= \begin{bmatrix}2(n_1-1)& 3m_1-2r_1& n_2& 2m_2\\3n_1-4& 4(m_1-r_1)& 2n_2& 3m_2\\n_1& 2m_1& 2(n_2-1)& 3m_2-2r_2\\2n_1& 3m_1& 3n_2-4& 4(m_2-r_2)\end{bmatrix}.\]
\end{itemize}

\end{theorem}
\begin{proof}
 By a suitable labelling of vertices of $G_1\ocirc G_2$, the distance matrix $D(G_1\ocirc G_2)$ can be written as
 \[\begin{bmatrix} 2(J_{n_1}-I_{n_1})& 3J_{n_1\times m_1}-2Q(G_1)& J_{n_1\times n_2}& 2J_{n_1\times m_2}\\3J_{m_1\times n_1}-2Q(G_1)^T& 4(J_{m_1}-I_{m_1})-2A(l(G_1))& 2J_{m_1\times n_2}& 3J_{m_1\times m_2}\\J_{n_2\times n_1}& 2J_{n_2\times m_1}& 2(J_{n_2}-I_{n_2})& 3J_{n_2\times m_2}-2Q(G_2)\\2J_{m_2\times n_1}& 3J_{m_2\times m_1}& 3J_{m_2\times n_2}-2Q(G_2)^T& 4(J_{m_2}-I_{m_2})-2A(l(G_2)) \end{bmatrix}\]\\
 
 Since $G_i$ is a $r_i$-regular graph $(i= 1, 2)$, it follows that $r_i$ is an adjacency eigenvalue of $G_i$ with an eigenvector $J_{{n_i}\times 1}$. If $U$ is a vector such that $A(G_1)U= \lambda_{1j}U$ (for $j= 2, 3, \ldots, n_1$), then $J_{{n_1}\times 1}^TU= 0$.\\

 Let $\phi= \begin{bmatrix} U\\Q(G_1)^TU\\0\\0\end{bmatrix}$ , then $D(G_1\ocirc G_2)\phi= -2(\lambda_{1j}+r_1+1) \phi$.\\
 
 Hence corresponding to each eigenvalue $\lambda_{1j}\neq r_1$ of $G_1$ we get an eigenvalue $-2(\lambda_{1j}+r_1+1)$ of $D(G_1\ocirc G_2)$ with an eigenvector $\phi.$\\
 
 Similarly, we can show that corresponding to each eigenvalue $\lambda_{2j}\neq r_2$ of $G_2$ there exists an eigenvalue $-2(\lambda_{2j}+r_2+1)$ of $D(G_1\ocirc G_2)$ with an eigenvector $\begin{bmatrix}0\\0\\W\\Q(G_2)^TW\end{bmatrix}$, where $W$ is the eigenvector corresponding to the eigenvalue $\lambda_{2j}\neq r_2.$\\
 
 Now let $Y$ be an eigenvector of $A(l(G_1))$ corresponding to the eigenvalue $\lambda_{1j}+r_1-2$ ($j= 2,\ldots, n_1$), then $J_{m_1\times 1}^TY=0.$
 If $\psi= \begin{bmatrix}Q(G_1)Y\\ -Y\\0\\0\end{bmatrix}$, then $D(G_1\ocirc G_2)\psi=0\psi$. This shows that $0$ is an eigenvalue of $D(G_1\ocirc G_2)$ with eigenvector $\psi.$\\
 
 Similarly, we can show that there exists an eigenvalue $0$ of $D(G_1\ocirc G_2)$ corresponding to the eigenvalue $\lambda_{2j}+r_1-2$ ($j=2,\ldots, n_2$).
 Let $Z$ be an eigenvector of $A(l(G_1))$, corresponding to the eigenvalue $-2$ of $A(l(G_1)).$\\
 
 If $\zeta=\begin{bmatrix}0\\Z\\0\\0\end{bmatrix}$, then $D(G_1\ocirc G_2)\chi= 0\zeta$. This implies that $0$ is an eigenvalue of $D(G_1\ocirc G_2)$ with multiplicity $m_1-n_1.$ 
 Similarly, corresponding to the eigenvalue $-2$ of $A(l(G_2))$, we get an eigenvalue $0$ of $D(G_1\ocirc G_2)$ with multiplicity $m_2-n_2$.\\
 Thus, in total we get $n_1+n_2+m_1+m_2-4$ eigenvalues of $D(G_1\ocirc G_2)$, then the remaining four eigenvalues are the eigenvalues of the equitable quotient matrix (Therefore, the total number of eigenvalues of $D(G_1\ocirc G_2)$ is  $n_1+n_2+m_1+m_2-4$, with the remaining four eigenvalues corresponding to the equitable quotient matrix)

 \[N=\begin{bmatrix}2(n_1-1)& 3m_1-2r_1& n_2& 2m_2\\3n_1-4& 4(m_1-r_1)& 2n_2& 3m_2\\n_1& 2m_1& 2(n_2-1)& 3m_2-2r_2\\2n_1& 3m_1& 3n_2-4& 4(m_2-r_2)\end{bmatrix}.\]
 
\end{proof}

 \begin{theorem}\label{se}
Let $G_i$ be a $r_i$-regular graph with $n_i$ vertices and $m_i$ edges. And let $A_i$ be the adjacency matrix of $G_i$ and  $\{r_i, \lambda_{i2},\ldots, \lambda_{in_{i}}\}$ ($i=1, 2$) be the adjacency spectrum of $G_i$. Then, the distance spectrum of $G_1\ominus G_2$ consists of the following,
\begin{itemize}
    \item [(i)]
$-\left(\lambda_{ij}+3\pm \sqrt{(\lambda_{ij}+1)^2+4(\lambda_{ij}+r_i)}\right), j=2,3,\ldots,n_i$,
    \item[(ii)] $-2$ with multiplicity $m_i-n_i$
    \item[(iii)] the four eigenvalues of the equitable quotient matrix of $D(G_1\ominus G_2)$,
\[N= \begin{bmatrix}4(n_1-1)-2r_1& 3m_1-2r_1& 3n_2& 2m_2\\3n_1-4& 2(m_1-1)& 2n_2& m_2\\3n_1& 2m_1& 4(n_2-1)-2r_2& 3m_2-2r_2\\2n_1& m_1& 3n_2-4& 2(m_2-1)\end{bmatrix}.\]
\end{itemize}

\end{theorem} 
\begin{proof}
    The proof follows by similar arguments in Theorem \ref{subdivision vertex}.
\end{proof}

(The next theorem provides new families of
$D$-cospectral graphs, the proof of which follows directly from Theorems \ref{subdivision vertex} and \ref{se}.)

\begin{theorem}
    (Let $G_1$ and $G_2$ be two $r$-regular $D$-cospectral graphs having the same number of vertices, and (let)$H$ be any arbitrary regular graph. )Then
\begin{itemize}
 \item[(i)] $ G_1\ocirc H$ and $ G_2\ocirc H$ are $D$-cospectral.
  \item[(ii)] $ G_1\ominus H$ and $ G_2\ominus H$ are $D$-cospectral.
\end{itemize}
\end{theorem}

The following theorem presents a new family of non $D$-cospectral $D$-equienergetic graphs of diameter $4$.

\begin{theorem}
  Let $G$ be a $r$-regular graph with $n$ vertices and $m$ edges. For a fixed $l\in \mathbb N $, let $\mathcal {P}_l$ be the family of all integer partitions on $'l'$ into parts, denoted by $p_1, p_2,\ldots, p_s$, each of size at least $3$. For $P= \{p_1, p_2,\ldots, p_s \}\in \mathcal {P}_l$, let $C_P$ be the union of cycles with vertices $p_1, p_2,\ldots, p_s$. Then the graphs $ G\ocirc C_P$, $P\in\mathcal{P}_l$ forms a family of $D$-equienergetic graphs.
\end{theorem}
\begin{proof}
Let $\rho_1=r\geq\rho_2\geq\dots\geq \rho_n$ and $\mu_1=2\geq \mu_2\geq\dots\geq\mu_l$ be the $A$-eigenvalues of $G$ and $C_P$ respectively. From Theorem \ref{subdivision vertex}, the $D$-eigenvalues of $ G\ocirc C_P$ are $-2(\rho_j+r+1)$, $j=2,3,\ldots,n$; $0$ with multiplicity $m-1$; $-2(\rho_{j_{x}}+3)$, $j_{x}=2,3,\ldots,l$, together with the eigenvalues of the equitable quotient matrix of $ G\ocirc C_P$,
 \[M= \begin{bmatrix}2(n-1)&3m-2r&l&2l\\3n-4&4(m-r)&2l&3l\\n&2m&2(l-1)&3l-4\\2n&3m&3l-4&4(l-2)\end{bmatrix}\] for every partition P of $'l'$.\\
 Since $\mu_{j_{x}}\geq -2$, for $j_{x}=2,3,\ldots,l$, we have $-2(\rho_{j_{x}}+3)\leq 0$.
Thus, $$\sum_{j_{x}=2}^l \lvert-2(\rho_{j_{x}}+3)\lvert= \sum_{j_{x}=2}^l 2(\rho_{j_{x}}+3)= 6l-10.$$\\
 Hence the energy remains the same for every partition of $'l'$.
\end{proof}
\begin{example}
    The graphs $K_{n,n}\ocirc C_6$ and  $K_{n,n}\ocirc C_3\cup C_3$ are of same number of vertices $3n+12$, with the following $D$-spectra,
    \begin{equation*}
        \begin{aligned}
            Spec_{D}(K_{n,n}\ocirc C_6)&=\{-2(n+1)^{(2(n-1))},-2,-4^{(4)},-10,0^{(n^2+4)},\alpha_1, \alpha_2,\alpha_3,\alpha_4\},\\
            Spec_{D}(K_{n,n}\ocirc C_3\cup C_3)&=\{-2(n+1)^{(2(n-1))},-2^{(2)},-4^{(2)},-8^{(2)},0^{(n^2+4)},\alpha_1, \alpha_2,\alpha_3,\alpha_4\},
        \end{aligned}
    \end{equation*}
    where $\alpha_1, \alpha_2,\alpha_3,\alpha_4$ are roots of the matrix \[\begin{bmatrix}4n-2& 3n^2-2n&6&12\\6n-4&4n^2-4n&12&18\\2n&2n^2&10&14\\4n&3n^2&14&16\end{bmatrix}.\]
    Clearly these two graphs are non $D$-cospectral $D$-equienergetic graphs of diameter $4$. Figure.$1$ illustrates the graphs when $n=1.$
\end{example}
\begin{figure}[htbp]
\centering
 \begin{tikzpicture}[scale=0.4,inner sep=1.3pt]
\draw (-2,0) node(3) [circle,draw,fill] {};
\draw (-2,8) node(1) [circle,draw,fill] {};
\draw (-2,4) node(2) [circle,draw,fill] {};
   \draw (4,7) node(5) [circle,draw,fill] {};
   \draw (4,11) node(4) [circle,draw,fill] {};
\draw (8,9) node(6) [circle,draw,fill] {};
\draw (4,9) node(10) [circle,draw] {};
\draw (6,8) node(11) [circle,draw] {};
\draw (6,10) node(12) [circle,draw] {};
\draw (4,-3) node(8) [circle,draw,fill] {};
\draw (4,-1) node(13) [circle,draw] {};
\draw (4,1) node(7) [circle,draw,fill] {};
\draw (8,-1) node (9) [circle,draw,fill] {};
\draw (6,0) node(15) [circle,draw] {};
\draw (6,-2) node(14) [circle,draw] {};
\draw (12,0) node(16) [circle,draw,fill] {};
\draw (12,8) node(17) [circle,draw,fill] {};
\draw (12,4) node(18) [circle,draw] {};
\draw (16,7) node(21) [circle,draw,fill] {};
\draw (16,4) node(22) [circle,draw] {};
\draw (16,1) node(23) [circle,draw,fill] {};
\draw (24,7) node(29) [circle,draw,fill] {};
\draw (24,4) node(28) [circle,draw] {};
\draw (24,1) node(27) [circle,draw,fill] {};
\draw (20,10) node(19) [circle,draw,fill] {};
\draw (20,-2) node(25) [circle,draw,fill] {};
\draw (18,8.5) node(20) [circle,draw] {};
\draw (18,-0.5) node(24) [circle,draw] {};
\draw (22,-0.5) node(26) [circle,draw] {};
\draw (22,8.5) node(30) [circle,draw] {};

\draw [-] (1) to (2) to (3) to (4) to (10) to (5) to (11) to (6) to (12) to (4) to (1) to (6) to (3) to (5) to (1);
\draw [-] (4) to (3) to (7) to (13) to (8) to (14) to (9) to (15) to (7) to (1) to (8) to (3) to (9) to (1);

\draw [-]  (16) to (18) to (17) to (19) to (20) to (21) to (22) to (23) to (24) to (25) to (26) to (27) to (28) to (29) to (30) to (19);
\draw[-](19) to (16) to (21);
\draw[-](23) to (16) to (25);
\draw[-](27) to (16) to (29);
\draw[-](21) to (17) to (23);
\draw[-](25) to (17) to (27);
\draw[-](29) to (17);

\end{tikzpicture}
\caption{$K_{1,1}\ocirc C_3\cup C_3,\; \text{and}\;K_{1,1}\ocirc C_6.$}
 \end{figure}

The next theorem determines new families of non $D$-cospectral $D$-equienergetic graphs of diameter $3.$
\begin{theorem}
    Let $G$ be a fixed triangle-free regular graph and $H$ be any regular graph with least eigenvalue at least $-2$. Then for any such graph $H$ with same order and regularity,
    \begin{itemize}
        \item [(i)]$G\dot{\vee} H$ forms a family of $D$-equienergetic graphs,
        \item [(ii)] $G\veebar H$ forms a family of $D$-equienergetic graphs.
    \end{itemize}
\end{theorem}
\begin{proof}
The proof is directly derived from Theorems $11$ and $12$ in \cite{central}.
\end{proof}

\section{Conclusion}
 In this paper, we investigated the distance spectral gap of connected graphs, providing various upper and lower bounds in terms of the sum of the squares of the $D$-eigenvalues. We also derived bounds for the distance eigenvalues and the distance energy of a graph $G.$ Furthermore, new families of non $D$-cospectral $D$-equienergetic graphs with diameters $3$ and $4$ are constructed. Additionally, new classes of non-isomorphic cospectral graphs are obtained.
\bibliographystyle{plain}
\bibliography{energy}

\end{document}